\documentclass[12pt]{article}
\usepackage{amssymb}
\usepackage{graphicx}
\usepackage{epsfig}

\def \QQ{Q \! \! \! \!\prime\;\,}
\def \CC{C \! \! \! \!\prime\;\,}
\begin{document}

\title{\bf Equivalent Binary Quadratic Form and the Extended Modular Group }

\author{M. Aslam Malik \thanks{aslam.math@pu.edu.pk} and  Muhammad Riaz \thanks{mriaz.math@pu.edu.pk} \\
Faculty of Science, Department of Mathematics\\
 University of the Punjab, Lahore, Pakistan.}

\date{}
\maketitle
\begin{abstract}

Extended modular group $\overline{\Pi}=\langle
R,T,U:R^2=T^2=U^3=(RT)^2=(RU)^2=1\rangle$,  where
 $ R:z\rightarrow -\overline{z}, ~T:z\rightarrow\frac{-1}{z}, ~U:z\rightarrow\frac{-1}{z +1}~$,
 has been used to study some properties of the binary quadratic forms whose base points lie in the point set
 fundamental region $F_{\overline{\Pi}}$ (See \cite{Tekcan1, Flath}).\\
  In this paper we look at how base points have been used in the study of equivalent binary quadratic forms,
  and we prove that two positive definite forms are equivalent if and only if the base point of one form is
  mapped onto the base point of the other form under the action of the extended modular group
  and any positive definite integral form can be transformed into  the reduced form of the same
  discriminant under the action of the extended modular group and extend these results for the
  subset $\QQ^*(\sqrt{-n})$ of the imaginary quadratic field $\QQ(\sqrt{-m})$. \\
\end{abstract}
AMS Mathematics subject classification (2000): 05C25, 11E04,11E25, 20G15\\
{\bf Keywords:} Binary quadratic forms, Extended modular group, imaginary quadratic field, Linear fractional transformations.\\
\section{Introduction}
The modular group $PSL(2,\mathbb {Z})$ has finite presentation
$\Pi=\langle T,U:T^2=U^3=1\rangle$  where
$T:z\rightarrow\frac{-1}{z}, U:z\rightarrow\frac{-1}{z +1}$ are
elliptic transformations and their fixed points in the upper half
plane are $i$ and $e^{2\pi i/3}$ \cite{Tekcan1}. The Modular group
$PSL(2,\mathbb {Z})$ is the free product of $C_2$ and $C_3$. If $n$
is finite, then the cyclic group of order $n$ is denoted by $C_n$
and an infinite cyclic group is denoted by $C_\infty$. The modular
group $PSL(2,\mathbb {Z})$  is the group of linear fractional
transformations of the upper half of the complex plane which have
the form $z\rightarrow\frac{az+b}{cz+d}$ where $a, b, c$ and  $d$
are integers, and $ad - bc = 1$, The group operation is given by the
composition of mappings.  This group of transformations is
isomorphic to the projective special linear group $PSL(2,\mathbb
{Z})$, which is the quotient of the 2-dimensional special linear
group over the integers by its center. In other words,
$PSL(2,\mathbb {Z})$  consists of all matrices
 \( \left(
\begin{array}{cc}
a & b \\
c & d
\end{array}
\right) \) where $a, b, c$ and $d$  are integers, and $ad - bc = 1$.\\
We assume the transformation $R(z)=-\overline{z}$ which represent the symmetry with respect to imaginary axis.
Then the group $\overline{\Pi}=\Pi\cup R\Pi$ is generated by the transformations $R,T$ and $U$ and has representation
$$\overline{\Pi}=\langle R,T,U:R^2=T^2=U^3=(RT)^2=(RU)^2=1\rangle. $$ $\overline{\Pi}$ is
called the extended modular group and $\Pi$ is a subgroup of index
$2$ in $\overline{\Pi}$. Therefore $\Pi$ is a normal subgroup of $\overline{\Pi}$ (See \cite{Tekcan1, Koruoglu}.\\
A binary quadratic form is a polynomial in two variable of the form
$$F=F(X,Y)=aX^2+bXY+cY^2$$
with real coefficients $a,b,c$. We denote $F$ simply by $[a,b,c]$. The discriminant of $F$ is denoted by
$\Delta(F)=b^2-4ac$. $F$ is an integral form if and only if $a,b,c \in \mathbb {Z}$, and $F$ is positive definite
 if and only if  $\Delta(F)<0$ and $a,c>0$. The form $F$ represent an integer $n$ if there are integers
  $x,y$  such that $F(x,y)=n$ and the representation of $n$ is primitive if $(x,y)=1$.\\
Lagrange was the first to introduce the theory of quadratic forms,
later Legendre expended the theory of quadratic forms, and greatly
magnified  even later by Gauss. It is proved by Gauss that the set
$C_\Delta$ of primitive reduced forms of discriminant $\Delta$
is an abelian group in a natural way.\\
Tekcan and Bizim   in \cite{Tekcan1}   discussed various properties
of the binary quadratic form in connection  with the extended
modular group.  Tekcan  in \cite {Tekcan2}  derived  cycles of the
indefinite quadratic forms and  cycles of the ideals. Continuing his
work on the quadratic forms, Tekcan  used base points of the
quadratic forms and derived  cycle and proper cycle of an indefinite
quadratic form  \cite {Tekcan3}. \\
 In \cite{Dani}  Dani and  Nogueira consider the actions of $SL(2, \mathbb{Z})$ and  $SL(2,
\mathbb{Z})_+$   on the projective space $\mathbb{P}$ and
 on  $\mathbb{P}\times \mathbb{P}$ and the results are obtained on orbits-closures to derive a class
  of binary quadratic forms.\\
The form $F=[a,b,c]$ is said to be almost reduced if $ a \leq c$ and
$\left|b\right| \leq a$ , with $a,c>0$ and $\Delta(F)<0$. This
condition is equivalent to $\left|Re(z)\right|\leq\frac{1}{2}$ and
$\left|z\right|\geq 1$. Thus $F$ is almost reduced if $\tau =z(F)$
lies in the fundamental region $F_\Pi$ described by
$\left|Re(z)\right|\leq\frac{1}{2}$ and $\left|z\right|\geq 1$ as shown in figure 1.\\

\begin{center}
\includegraphics[width=3.5in, height=3.0in]{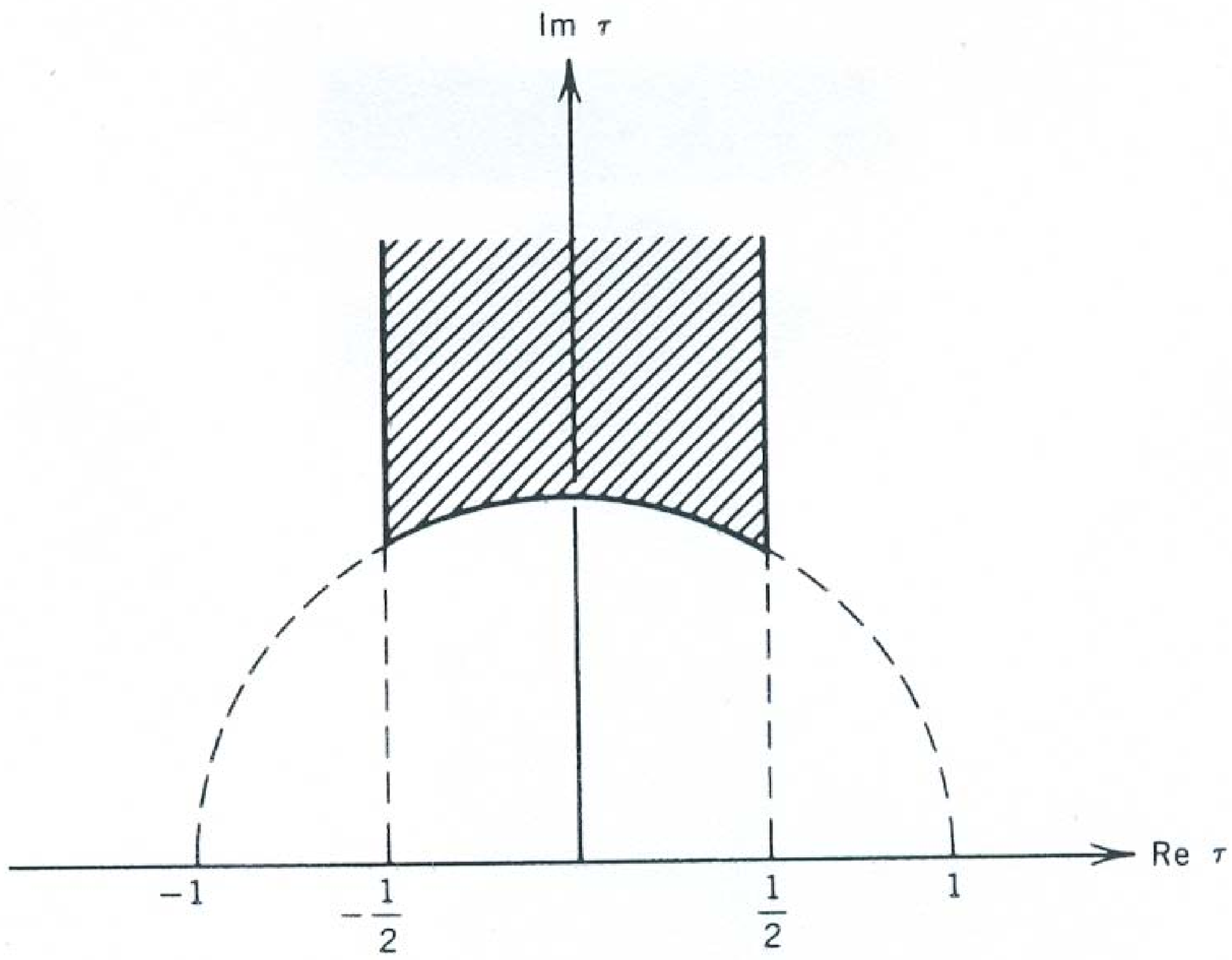}
\end{center}
\centerline{Fig: 1}

The form $F=[a,b,c]$ is said to be reduced if $\left|b\right|\leq a \leq c$, and in case $a=\left|b\right|$, then $a=b$
and in case $a=c$ then  $b\geq 0$ \cite{Flath}, then the corresponding fundamental region is denoted
by $F_{\overline{\Pi}}$ as shown in figure 2.\\

\begin{center}
\includegraphics[width=3.5in, height=3.0in]{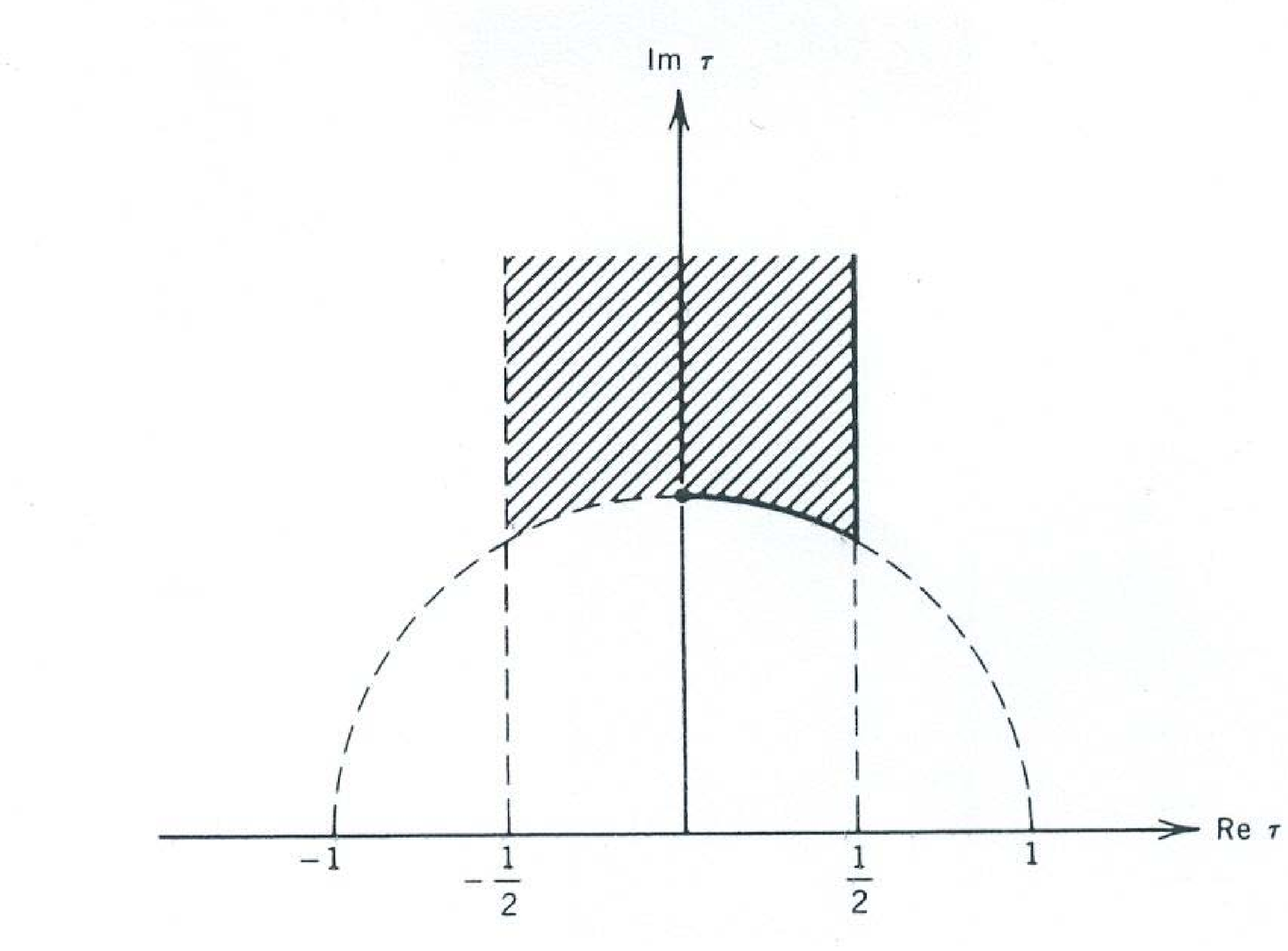}
\end{center}
\centerline{Fig: 2}

Flath \cite{Flath} proved that the number of almost reduced forms
with discriminant $\Delta(F)=b^2-4ac$ is finite. In \cite{Imrana} it
has been proved that $\alpha$ is mapped onto $\overline{\alpha}$
under the action of $PSL(2,\mathbb {Z})$  if and only if the
quadratic form $f$ is equivalent to $-f$. It is well known that the
set
$$\QQ^*(\sqrt{n}):=\{\frac{a+\sqrt{n}}{c}:a,c,b=\frac{a^2-n}{c}\in
\mathbb{Z}~\textmd{and}~(a,b,c)=1\}$$
is a $G$-subset of the real quadratic field  $\QQ(\sqrt{m})$ under the action of $PSL(2,\mathbb {Z})$.
 In \cite{Asim} it was proved that there exist two proper $G$-subsets of $\QQ^*(\sqrt{n})$
 when $n\equiv0~(mod~p)$ and four $G$-subsets
of $\QQ^*(\sqrt{n})$ when $n\equiv0~(mod~pq)$. \\
 Malik and Riaz in \cite{Riaz1}  extended this idea to determine four
proper $G$-subsets of $\QQ^*(\sqrt{n})$ with $n\equiv0~(mod~2pq)$,
We generalized this result for $n\equiv 0~(mod~p_1 p_2...p_r)$  and
proved that there are exactly $2^r$ number of $G$-subsets of
$\QQ^*(\sqrt{n})$. We also proved for $h=2k+1\geq3$  there are
exactly two $G$-orbits of $\QQ^*(\sqrt{2^h})$ namely $
(2^k\sqrt{2})^G$ and $(\frac{2^k\sqrt{2}}{-1})^G$. In the same paper
we used subgroup $G^*=\langle yx\rangle$ and $G^{**}=\langle yx,y^2x
\rangle$ to determine the $G$-subsets and $G$-orbits of $\QQ^*(\sqrt{n})$.\\
If $p$ is an odd prime, then  $\lambda~{\not\equiv}~0(mod~p)$ is said to
be a quadratic residue of $p$ if there exists an integer $r$ such that $r^2\equiv \lambda (mod~p)$.\\
The quadratic residues of $p$ form a subgroup $R$ of the group of nonzero integers modulo $p$ under multiplication and
$|R|=(p-1)/2$.  \cite{Adler}\\
The norm of an element $\alpha$ in $Q^*(\sqrt{n})$ is defined by
$N(\alpha)=\alpha\overline{\alpha}= \frac{a^2-n}{c^2}$.\\
The Legendre symbol $\left(a/p\right)$ is defined as $1$ if $a$ is a quadratic residue of
$p$ otherwise it is defined by $-1$.\\
\\
\textbf{Theorem 1.1} \cite{Adler}\\
(i) Let $p$ be any odd prime such $p\equiv1~(mod~4)$ and $a$ be quadratic residue of $p$ then $p-a$ is a quadratic residue of $p$.\\
(ii) Let $p$ be any odd prime such $p\equiv3~(mod~4)$ and $a$ be quadratic residue of $p$ then $p-a$ is a quadratic non-residue of $p$.\\
\textbf{Theorem 1.2} \cite{Flath}\\
The number of almost reduced forms $\underline{h}\left( \Delta \right)$ , with given discriminant $ \Delta <0$, is finite.\\

\section{ The relation between quadratic forms and the extended modular group.}
It is interesting to determine equivalent binary quadratic forms,
the interest has risen in the last decade as the possibilities to
study binary quadratic forms in combination with  mathematical
software has improved. There is a strong relation between equivalent
binary quadratic forms and their base points  (See \cite {Tekcan1, Tekcan3}).\\
Two forms $F=[a,b,c]$ and $G=[A,B,C]$ are said to be equivalent if
there exists an element g=\( \left(
\begin{array}{cc}
r & s \\
t & u
\end{array}
\right) \) $\in \bar{\Pi}$, where $r,s,t,u$ are integers and
$ru-st=1$, such that
 $gF(X,Y)=G(X,Y)$\\
or $F(rX+sy,tX+uY)=G(X,Y)$\\
the coefficients of $G$ in terms of the coefficients of $F$ are as follows\\
$A=ar^2+brt+ct^2$\\
$B=2ars+bru+bst+2ctu$\\
$C=as^2+bsu+cu^2.$\\
Then it is easy to see that $\Delta(F)=(ru-st)^2 \Delta(gF)$ but
$ru-st=1$ implies that
 $\Delta(F)=\Delta(gF)$ \cite {Mollin}, that is any two equivalent forms have the same discriminant,
 but the converse is not true in general.
The definition of $gF$ is a group action on the set of forms because it satisfies the two axioms\\
(i) \( \left(
\begin{array}{cc}
1 & 0 \\
0 & 1
\end{array}
\right) \) $F=F$  $~~~$ (ii) $g(hF)=(gh)F.$\\
Gauss (1777-1855) introduced the group action of extended modular
group $ \bar{\Pi}$ on the set of forms.  We use the standard
notations as used in \cite{Tekcan1}. Let $\mathbb {U}=\{z\in
\CC:Im(z)>0~\}$ be the upper half of the complex plane. For a
positive definite form $F=[a,b,c]$ with $\Delta(F)<0$ there is a
unique $z=x+iy=z(F) \in \mathbb {U}$ such that
$$F=a(X+zY)(X+\overline{z}Y)$$
 then $z=x+iy$ is called base point of $F$ in $\mathbb {U}$.\\
 $$F=a(X+zY)(X+\overline{z}Y)=aX^2+2axXY+a\left|z\right|^2Y^2$$
By comparison we get  $x=\frac {b}{2a}$ and $y=\frac{1}{2a}
\sqrt{-\Delta(F)}>0$
so $z=x+iy=\frac{b+i\sqrt{-\Delta(F)}}{2a}\in \mathbb {U}$ .\\
Conversely for any $w=x+iy\in \mathbb {U}$ there exists a positive
definite quadratic form $F=\left[\frac{1}{\left|w\right|^2},
\frac{2x}{\left|w\right|^2}, 1\right]$ with base point $w$ and
discriminant $\Delta(F)<0$.\\
 There is a one-one correspondence
between the set of positive definite forms with fixed
 discriminant and the set of base points of these forms  \cite{Tekcan1}.\\
 In the next theorem we show that there is a
strong relation between equivalent binary quadratic forms and their
base points. We proved that two positive definite
 forms are equivalent if and only if the base point of one form is mapped onto the base point of the other
  form under the action of the extended modular group.\\
\textbf{Theorem 2.1}\\
Two positive definite forms $F=[a,b,c]$ and $G=[A,B,C]$ are equivalent if and only if
$ \alpha= g(\beta), g \in \overline{\Pi}$, where  $\alpha$ and $\beta$ are the base points of $F$ and $G$ respectively. \\
\textbf{Proof.}\\
Let $F=[a,b,c]=aX^2+bXY+cY^2$ and $G=[A,B,C]=Au^2+Buv+Cv^2$ are  equivalent and $\alpha=z(F)$ and $\beta=z(G)$ are the base points of
$F$ and $G$ respectively, then we can assume that\\
$$F=a(X+\alpha Y)(X+\overline{\alpha}Y)$$
$$G=A(u+\beta v)(u+\overline{\beta}v)$$
Since $F$ and $G$ are equivalent then exist an element
g=\( \left(
\begin{array}{cc}
r & s \\
t & u
\end{array}
\right) \) $\in \overline{\Pi}$, where $r,s,t,u$ are integers and
$ru-st=1$, such that $u=rX+sY$ and $v=tX+uY$.
Thus by \cite{Mushtaq} we have\\
$a(X+\alpha Y)(X+\overline{\alpha}Y)=A(u+\beta v)(u+\overline{\beta}v)$\\
$=A(rX+sY+\beta (tX+uY))(rX+sY+\overline{\beta} (tX+uY))$\\
$=A(rX+sY+\beta tX+ \beta uY)(rX+sY+\overline{\beta} tX+ \overline{\beta}uY)$\\
$=A((r+\beta t)X+(s+\beta u)Y)((r+\overline{\beta} t)X+(s+\overline{\beta} u)Y)$\\
$=A'(X+\left(\frac{s+\beta u}{r+\beta t}\right) Y)(X+\left(\frac{s+\overline{\beta} u}{r+\overline{\beta}t}\right) Y)$\\
where $A'=A(r+\beta t)(r+\overline{\beta} t)$.\\
Thus, if $\alpha=\frac{s+\beta u}{r+\beta t}$, then $\overline{\alpha}=\frac{s+\overline{\beta} u}{r+\overline{\beta}t}$\\
or, if  $\alpha=\frac{s+\overline{\beta} u}{r+\overline{\beta}t}$, then $\overline{\alpha}=\frac{s+\beta u}{r+\beta t}$
It is clear that $\alpha =g'(\beta)$ and  $\overline{\alpha} =g'(\overline{\beta})$ under the linear fractional transformation
$$g'(z)=\frac{s+z u}{r+zt},~\textmd{where}~ st-ru=1$$
Conversely suppose that there exists a linear fractional transformation
 $g(z)=\frac{s+z u}{r+zt} \in \overline{\Pi},~\textmd{where}~ st-ru=1$
 such that $\alpha= g(\beta)$, then $\alpha$ and $\overline{\alpha}$ are images of
   $\beta$ and $\overline{\beta}$ under $g$, where $\alpha$ and $\beta$ are the base
    points of $F$ and $G$ respectively, then clearly $ gF=G$. This proves the result. \quad\quad$\Box$ \\
The following results are the immediate consequences of the above theorem.\\
\textbf{Corollary 2.2}\\
Two forms $F$ and $G$ are equivalent if and only if their base points lie in the same fundamental region $F_{\overline{\Pi}}$.\\
\textbf{Corollary 2.3}\\
If $F$ and $G$ are two equivalent binary quadratic forms then  both forms are either almost
 reduced quadratic forms or reduced quadratic forms.\\
\textbf{Theorem 2.4}\\
The binary quadratic form $F=x^2+py^2$ is equivalent to $\lambda F=\lambda x^2+ \lambda py^2$
 if and only if $\left(\lambda/ p\right)=1$, For any $\lambda \in \mathbb{Z}$ and $p$ be any prime.\\
\textbf{Proof.}\\
Let the binary quadratic form $F$ be equivalent to  $\lambda F$, then there exists an element
g=\( \left(
\begin{array}{cc}
r & s \\
t & u
\end{array}
\right) \) $\in \bar{\Pi}$, where $r,s,t,u$ are integers and
$ru-st=1$, such that
 $gF(x,y)= \lambda F(x,Y)$\\
 $\Leftrightarrow F(rx+sy,tx+uY)= \lambda F(x,y)$\\
 $ \Leftrightarrow (rx+sy)^2+p(tx+uy)^2= \lambda x^2+ \lambda py^2$ (Comparing the coefficients of $x^2$)\\
$ \Leftrightarrow  r^2+pt^2=\lambda $\\
$ \Leftrightarrow r^2 \equiv \lambda~(mod~p)$.\\
Thus  $\lambda$ is a quadratic residue modulo $p$.  \quad\quad$\Box$\\
\textbf{Example 2.5}\\
The form $x^2+37y^2$ is equivalent to the form $-x^2-37y^2$ because $\left(\lambda/ 37\right)=1$
but the form $x^2+79y^2$ is not  equivalent to the form $-x^2-79y^2$ because $\left(\lambda/ 79\right)=-1$. \quad\quad$\Box$\\
\textbf{Theorem 2.6}\\
Let $F$ be a positive definite integral form with $\alpha=z(F)$ as the base points of $F$,
and let $g(\alpha)=\beta$ under the action of the extended modular group. Then $F$ is equivalent
 to its reduced form $F_R$ whose base point is $\beta=z(F_R)$.\\
\textbf{Proof.}\\
Let $F=ax^2+bxy+cy^2$ be primitive positive definite form and let $n$ be the smallest integer
represented by $F$, put  $F_R=Au^2+Buv+Cv^2$ and  $\alpha=z(F)$ and $\beta=z(F_R)$ be the base
 points of $F$ and $F_R$ respectively, If $g(\alpha)=\beta$ then
\( \left(
\begin{array}{cc}
r & s \\
t & u
\end{array}
\right) \) $(\alpha)=\beta \Rightarrow$
$\frac{r \alpha  +s}{ t \alpha  +u}=\beta$\\
$r \alpha  +s = \beta \left( t \alpha  +u \right)$\\
$\Rightarrow r \alpha  +s = \beta u+ t\alpha\beta  \Rightarrow  r \alpha  +s - \beta u- t\alpha\beta =0 $\\
$\Rightarrow \left( r \alpha  +s \right)-  \left(u+ t\alpha \right)\beta=0$\\
$ r \alpha +s =0$ and $\left(u+ t\alpha \right)\beta=0$,  $\beta\neq 0$\\
$ r \alpha +s =0$ and $u+ \alpha t=0$\\
$\frac{s}{-\alpha r} =\frac{u}{-\alpha t} \Rightarrow s=-\alpha r$ and $u=-\alpha t$\\
$\frac{s}{u} =\frac{-\alpha r}{-\alpha t} \Rightarrow \frac{s}{u} =\frac{r}{t} $\\
$ru=st \Rightarrow ru-st=1 $ this shows that $g \in \bar{\Pi}$, put $u=rx+sy$ and $v=tx+uy$
then $F$ is transformed into the form $F_R$ as  follows:\\
$F_R(x,y)= (ar^2+brt+ct^2)x^2 +(2ars+bru+bst+2ctu)xy+(as^2+bsu+cu^2)y^2$\\
or $F_R(x,y)= nx^2+Bxy+Cy^2$, where\\
 $n=ar^2+brt+ct^2$\\
  $B=2ars+bru+bst+2ctu$ \\
   $C=as^2+bsu+cu^2$.\\
Then $\Delta(F)=(ru-st)^2 \Delta(F_R)$ this implies that
$\Delta(F)=\Delta(F_R)$. Now by definition $a\leq c \Rightarrow 4a^2
\leq 4ac = b^2- \Delta(F) \leq a^2- \Delta(F)$, where $\Delta(F) <0$
  implies that $ 3a^2\leq \left(-\Delta(F) \right)$ and finally $a\leq \sqrt{\left(-\Delta(F)\right)/3}$,
   Thus by \cite{Mollin}, we see that $\left|B\right| <n$ then  $F$ is properly equivalent to $F_R$,
    since $F$ is positive definite, and $F_R(0,1)=C$ which implies that $C\in \mathbb{N}$,
     and $C\geq n$ by the minimality of $n$.  This proves the result.  \quad\quad$\Box$\\
\section{ The relation between quadratic forms and orbits of imaginary quadratic fields.}
The imaginary quadratic fields are defined by the set $\QQ(\sqrt{-m})=\{a+b\sqrt{-m}:a,b\in \QQ \}$, where $m$ is a
square free positive integer. We shall denote the subset
 $$\{\frac{a+\sqrt{-n}}{c}:a,c,b=\frac{a^2+n}{c}\in
\mathbb {Z}, c\neq 0 \}$$
by $\QQ^*(\sqrt{-n})$. The imaginary quadratic fields are very useful in different branches of Mathematics.
 The integers in  $\QQ(\sqrt{-1})$  are called Gaussian integers, and the integers in $\QQ(\sqrt{-3})$
  are called  Eisenstein integers.\\
We denote the modular group $\Pi$ by $G$ for our convenience and use
coset diagrams to investigate  $G$-subsets and $G$-orbits of $\QQ^*(\sqrt{n})$ under the action of the modular group $\Pi$ (See \cite{Asim, Riaz1}).\\
   In \cite{Imrana} it has been proved that for any quadratic form $F=Ax^2+Bxy+Cy^2$ if $\alpha$ and
   $ \overline{\alpha}$ belong to   the subset $\QQ^*(\sqrt{n})$  of the real quadratic field then there exist a
   rational number $\lambda$ such that $F(x,y)=\lambda (cx^2-2axy+by^2)$. We extend this result for the subset
    $\QQ^*(\sqrt{-n})$ of the imaginary quadratic field in the next Lemma.\\
    \\
\textbf{Lemma 3.1}\\
For any positive definite binary quadratic form $F=[A,B,C]=Ax^2+Bxy+Cy^2$ if
$\alpha=\frac{a+\sqrt{-n}}{c} \in \QQ^*(\sqrt{-n})$ with $ b=\frac{a^2+n}{c}$
 such that $\alpha$ be the base point of $F$ .
 Then there exists a rational number $A'$ such that $F(x,y)=A'(cx^2-2axy+by^2)$.\\
\textbf{Proof.}\\
Let $\alpha=\frac{a+\sqrt{-n}}{c} \in \QQ^*(\sqrt{-n})$ with $ b=\frac{a^2+n}{c}$ such that $\alpha$ is base point of $F$ then
$F=A(x+\alpha y)(x+\overline{\alpha}y)$ or $F=A(x+(\frac{a+\sqrt{-n}}{c}) y)(x+(\frac{a-\sqrt{-n}}{c})y)$
 whose roots are $~- \frac{a\pm \sqrt{-n}}{c}$ and  again the roots of $F=[A,B,C]=Ax^2+Bxy+Cy^2$ are
 $ \frac{-B\pm \sqrt{B^2-4AC}}{2A}$. By comparing these roots the result
  follows immediately by \cite{Imrana}. \quad\quad$\Box$\\
Further  it has been proved in \cite{Imrana} that $\alpha$ is mapped
onto
 $\overline{\alpha}$ under the action of $PSL(2,\mathbb {Z})$ on $\QQ^*(\sqrt{p})$,
  where $\alpha=\frac{a+\sqrt{p}}{c},~ b=\frac{a^2-p}{c}~\textmd{and}~(a,b,c)=1$
  if and only if the quadratic form $f(x,y)=cx^2-2axy+by^2$ is equivalent to $-f$.
  We generalized this result for the subset $\QQ^*(\sqrt{-n})$ of the imaginary quadratic field.\\
\textbf{Theorem 3.2}\\
Under the action of $PSL(2,\mathbb {Z})$ on $\QQ^*(\sqrt{-n})$, $\alpha$ is mapped onto $\beta$
if and only if the  binary quadratic forms $F$ and $G$ are equivalent, where $\alpha=z(F)$ and $\beta=z(G)$ are the base points of
$F$ and $G$ respectively.\\
\textbf{Proof.}\\
We know that the modular group $PSL(2,\mathbb {Z})$ is a subgroup of index $2$ in $\overline{\Pi}$,
also $\Pi$ is a normal subgroup of $\overline{\Pi}$, Let $\alpha=\frac{a+\sqrt{-n}}{c} \in \QQ^*(\sqrt{-n})$
 with $ b=\frac{a^2+n}{c}$ such that $\alpha$ be the base point of $F$ and $\beta$ be the base points of $G$ respectively.
  Since $\alpha=g(\beta)$ under the action of  $PSL(2,\mathbb {Z})$ on $\QQ^*(\sqrt{-n})$, then by Theorem 2.1, $F$ is equivalent to $G$.
  This proves the result.  \quad\quad$\Box$\\
\textbf{Conclusion}\\
 There is a one-one correspondence between the set of positive definite forms with fixed
 discriminant and the set of base points of these forms.
We prove that there is a strong relation between equivalent binary
quadratic forms and their base points in a sense that two positive
definite  forms are equivalent if and only if the base point of one
form is mapped onto the base point of the other   form under the
action of the extended modular group and two forms are equivalent if
and only if  their base points lie in the same fundamental region.
These results can be extended to determine the quadratic ideals for
a given quadratic irrationals, and for determining proper cycles of
an indefinite quadratic forms as done in \cite{Tekcan2, Tekcan3}.  \\

\end{document}